\documentclass[11pt]{amsart}
\usepackage{amsmath,amssymb,amsthm}
\usepackage{needspace}
\usepackage[hidelinks]{hyperref}
\newtheorem{theorem}{Theorem}[section]
\newtheorem{lemma}[theorem]{Lemma}
\newtheorem{proposition}[theorem]{Proposition}
\newtheorem*{lewkotheorem}{Theorem (Lewko)}
\newcommand{\F}{\mathbb F}
\newcommand{\PP}{\mathbb P^1}
\newcommand{\PGL}{\operatorname{PGL}_2}
\newcommand{\card}[1]{\lvert #1\rvert}
\newcommand{\one}{\mathbf 1}
\newcommand{\op}{\circ}
\title[M\"obius expansion and pointwise energy transfer]{M\"obius expansion and pointwise energy transfer over arbitrary fields}
\author{Yao Zhi}
\email{ericyao2026@gmail.com}
\date{}
\begin{document}
\begin{abstract}
Let $F$ be a field of positive characteristic $p$, and let
$A\subset F^\times$ have $N$ elements. We prove that, for every
$\lambda\in F^\times$,
\[
 \card{A+\lambda A^{-1}}\gg\min\{N^{10/9},(pN^2)^{1/3}\}.
\]
This improves the exponent $31/30$ to $10/9$ in the common small-set
range. Over prime fields the bound strengthens to
$\min\{N^{10/9},(pN)^{1/2}\}$. Analogous bounds hold for $Af(A)$ when $f$ is a M\"obius
transformation satisfying the stated nondegeneracy conditions.
In characteristic zero, the bound $N^{10/9}$ remains valid.
No subfield-avoidance condition is required.
\end{abstract}
\maketitle

\section{Introduction}
A basic problem in additive combinatorics is to estimate the
growth of sets formed by combining a finite set with its image
under a nonlinear map. For the reciprocal map, this leads to
\[
 A+\lambda A^{-1}
 =\{a+\lambda b^{-1}:a,b\in A\},\qquad
 A\subset F^\times,\quad \lambda\in F^\times.
\]
Aksoy Yazici, Murphy, Rudnev and Shkredov \cite{AYMRS} proved
$\card{A+A^{-1}}\gg N^{31/30}$, where $N=\card A$, for sets in a
positive-characteristic small-set range. For $A\subset\F_p^\times$,
Balog, Broughan and Shparlinski
\cite[Theorem~7 and Corollary~9]{BBS} proved
\[
 \card{A+A^{-1}}\gg
 \min\left\{\sqrt{pN},\frac{N^2}{\sqrt p}\right\}.
\]
Their three-set estimate also applies to $A+\lambda A^{-1}$
for every $\lambda\ne0$.

Write $p=\operatorname{char}F$ in positive characteristic, and
$p=\infty$ in characteristic zero. For a M\"obius transformation put
\[
 f(x)=\frac{\alpha x+\beta}{\gamma x+\delta},\qquad
 \alpha\delta-\beta\gamma\ne0.
\]
\Needspace{14\baselineskip}
Our main result is the following.

\begin{theorem}\label{thm:main}
Let $F$ be a field, and let $A\subset F^\times$ be a finite
nonempty set with $\card A=N$. Let $S$ be one of the following sets:
\begin{enumerate}
\item[(i)] $S=A+\lambda A^{-1}$, where $\lambda\in F^\times$;
\item[(ii)] $S=Af(A)$, where $\alpha\beta\gamma\delta\ne0$
 and $f$ is defined and nonzero on $A$.
\end{enumerate}
In positive characteristic $p$ we have
\begin{equation}\label{eq:main}
 \card S\gg \min\{N^{10/9},(pN^2)^{1/3}\}.
\end{equation}
For $F=\F_p$ the stronger estimate holds:
\begin{equation}\label{eq:prime}
 \card S\gg \min\{N^{10/9},(pN)^{1/2}\}.
\end{equation}
In characteristic zero, $\card S\gg N^{10/9}$. All implied
constants are absolute.
\end{theorem}

In the common range $p>2$ and $N<p^{5/8}$, the exponent $10/9$
improves the exponent $31/30$ of Aksoy Yazici, Murphy, Rudnev
and Shkredov \cite{AYMRS}. Over $\F_p$, the present bound is
stronger than \cite{BBS} for $N<p^{9/16}$, weaker for
$p^{9/16}<N<p^{9/11}$, and of the same $\sqrt{pN}$ scale for
$N\ge p^{9/11}$. Over the reals, Stevens and Warren
\cite[Corollary~7]{SW} obtained the larger exponent $49/38$,
up to logarithmic factors. For sufficiently large subsets of prime
fields, Semchankau \cite{Sem} obtained sharper asymptotic estimates.

The point--line incidence estimate used here is due to Lewko
\cite{Lewko}. Independently, Miao and Xie \cite{MX} obtained a
related estimate over prime fields, and Yu \cite{Yu} subsequently
studied refinements of the polynomial-method approach.

Related composition-energy estimates for M\"obius transformations
appear in the work of Shkredov \cite[Lemma~14]{Shkredov},
Rudnev and Wheeler \cite{RW}, and Warren and Wheeler \cite{WW}.

We briefly outline the argument. Let $\Gamma$ be either $(F,+)$
or $(F^\times,\cdot)$, write its operation as $\op$, and denote
the inverse of $x$ by $x^{-\op}$. For a finite set
$Y\subset\Gamma$, put
\[
 r_Y(t)=\card{\{(y_1,y_2)\in Y^2:
             y_1\op y_2^{-\op}=t\}},\qquad
 E_2(Y)=\sum_t r_Y(t)^2.
\]
Given finite sets $X,Y\subset\Gamma$, consider the family of
M\"obius transformations
\[
 h_{a,z}(x)=a\op f(z\op x),\qquad
 a\in X,\quad z\in Y^{-\op}.
\]
Fixing an incident pair and changing projective coordinates
reduces the corresponding M\"obius transformations to affine maps,
so Lewko's theorem bounds their rich incidences. The relative
compositions split into exceptional and regular parameter fibres,
with the multiplicities of the latter controlled by $E_2(Y)$.

\Needspace{9\baselineskip}
Each representation $x=y_1\op y_2$ with $y_1,y_2\in Y$
contributes $\card X$ incidences at $x$. Comparing this pointwise
count with the weighted second moment gives
Proposition~\ref{prop:asymmetric}, which controls $E_2(Y)$
through $\card{X\op f(Y)}$ without assuming equal input sizes.
Applying it to both inputs and using Cauchy--Schwarz proves
Theorem~\ref{thm:main}. Over $\F_p$, centering on
$\mathbb P^1(\F_p)$ gives the stronger bound \eqref{eq:prime}.

The paper is organized as follows. Section~2 derives estimates for
rich M\"obius transformations, including their centered counterparts
over prime fields. Section~3 bounds the composition multiplicities of
the relevant transformation families and obtains a weighted incidence
estimate. Section~4 establishes the pointwise energy-transfer estimate
and proves Theorem~\ref{thm:main}.

\section{Rich M\"obius transformations and centered estimates}
We recall the point--line incidence estimate of Lewko
\cite[Theorem~1]{Lewko} in the form needed below.

\begin{lewkotheorem}
Let $F$ be a field of characteristic $p$, with $p=\infty$ in
characteristic zero. Let $P\subset F^2$ be a finite set of points
and let $\mathcal L$ be a finite family of distinct affine lines.
Then
\begin{equation}\label{eq:ST}
 I(P,\mathcal L)\le 3(\card P\card{\mathcal L})^{2/3}
 +\card P+\card{\mathcal L}+2\card P\card{\mathcal L}/p.
\end{equation}
If $F=\F_p$, the coefficient $2$ in the final term can be
replaced by $1$.
\end{lewkotheorem}

For a finite set $S\subset F$ with $\card S=M$, let
\[
 i(g)=\card{\{x\in S:g(x)\in S\}}\qquad(g\in\PGL(F)).
\]
All transformations are interpreted on the projective line; in particular,
a pole is not counted as an incidence with $S$.

\begin{lemma}\label{lem:rich}
For $k\ge C(1+M^2/p)$,
\begin{equation}\label{eq:tail}
 \card{\{g:i(g)\ge k\}}\ll M^6/k^4.
\end{equation}
Here distinct transformations are counted only once. Moreover, for every finite family $\mathcal G$,
\begin{equation}\label{eq:richsum}
 \sum_{g\in\mathcal G}i(g)\ll M^{3/2}\card{\mathcal G}^{3/4}+(1+M^2/p)\card{\mathcal G}.
\end{equation}
If $F=\F_p$, set $q=p+1$ and $\mu=M^2/q$. Then, for $k\ge C$,
\begin{equation}\label{eq:centertail}
 \card{\{g:i(g)-\mu\ge k\}}\ll M^6/k^4,
\end{equation}
and
\begin{equation}\label{eq:centersum}
 \sum_{g\in\mathcal G}(i(g)-\mu)_+\ll M^{3/2}\card{\mathcal G}^{3/4}+\card{\mathcal G}.
\end{equation}
\end{lemma}
\begin{proof}
It suffices to consider $M\ge2$: when $M=1$, the claims follow from $i(g)\le1$ and the linear terms in the bounds. The large-deviation sets are empty for $k>M$, so we further restrict to $k\le M$.

Fix an incident pair $(a,b)\in S^2$ and restrict attention to transformations $g$ satisfying $g(a)=b$. Apply the coordinate changes $x\mapsto 1/(x-a)$ on the source and $y\mapsto 1/(y-b)$ on the target. The resulting transformations fix the point at infinity and are therefore affine. Distinct original transformations give distinct affine maps, hence distinct nonvertical lines. All other incidences lie in a Cartesian grid of $(M-1)^2$ points. If $h$ such transformations each have at least $k$ incidences before conjugation, \eqref{eq:ST} yields
\[
 (k-1)h\le 3((M-1)^2h)^{2/3}+(M-1)^2+h
              +\frac{2(M-1)^2h}{p}.
\]
For $k\ge C(1+M^2/p)$ with a sufficiently large absolute constant $C$,
the terms linear in $h$ can be absorbed, giving
\[
 kh\ll M^{4/3}h^{2/3}+M^2,
 \qquad h\ll \frac{M^4}{k^3}+\frac{M^2}{k}
 \ll\frac{M^4}{k^3},
\]
The last inequality uses $k\le M$. Each map counted in \eqref{eq:tail} belongs to at least $k$ of the $M^2$ fixed-pair families. Summing over these families and dividing by $k$ proves \eqref{eq:tail}. Over an infinite field, the set of maps under consideration is finite for $k\ge3$, because the images of three distinct points determine a projective linear transformation.

For the prime-field assertion put $q=p+1$, $\mu=M^2/q$ and
$\mu_0=(M-1)^2/p$. Since $M\le p$, the exact identity
\[
 \mu-1-\mu_0=-\frac{(p+1-M)^2}{p(p+1)}\in[-1,0]
\]
shows that a transformation with $i(g)-\mu\ge k$ and $g(a)=b$ has at least $\mu_0+k-1$ incidences in the remaining grid. The final term of \eqref{eq:ST} has coefficient exactly $1$ over prime fields, so the corresponding $h$ lines satisfy
\[
 (\mu_0+k-1)h
 \le 3((M-1)^2h)^{2/3}+(M-1)^2+h+\mu_0h.
\]
After cancelling $\mu_0h$, the residual $O(h)$ term can be absorbed when $k\ge C$, giving $h\ll M^4/k^3$ as before. Since $i(g)\ge\mu+k\ge k$, each such transformation belongs to at least $k$ fixed-pair families, which proves \eqref{eq:centertail}. This cancellation requires the exact coefficient $1$ and does not follow from an estimate with coefficient $2$.

Finally, for any finite $\mathcal G$ the layer-cake formula and
\eqref{eq:tail} yield, with $K\asymp 1+M^2/p$,
\[
 \sum_{g\in\mathcal G}i(g)
 \ll K\card{\mathcal G}
  +\sum_{k\ge K}\min\left\{\card{\mathcal G},\frac{M^6}{k^4}\right\}.
\]
Split the sum at $k_0=(M^6/\card{\mathcal G})^{1/4}$. If $k_0<K$, the tail sum is absorbed by $K\card{\mathcal G}$. This proves \eqref{eq:richsum}. Applying the same argument to $(i(g)-\mu)_+$ with $K\asymp1$ establishes \eqref{eq:centersum}. The case of an empty family is immediate.
\end{proof}

\section{Composition multiplicities}
Let $U,Q\subset\Gamma$ be nonempty finite sets of sizes $n,u$. Consider the family
\begin{equation}\label{eq:h}
 \mathcal H=\{h_{a,z}:a\in U,z\in Q\},\qquad h_{a,z}(x)=a\op f(z\op x).
\end{equation}
The group operations in this definition extend to translations or dilations of $\PP(F)$. Define the associated composition multiplicities by
\[
 w(g)=\card{\{(h_1,h_2)\in\mathcal H^2:h_1h_2^{-1}=g\}}.
\]
\begin{lemma}\label{lem:weights}
The transformations in \eqref{eq:h} are pairwise distinct. Their composition weights admit a decomposition into exceptional and regular contributions, defined at the level of parameter tuples, such that
\begin{align}
 \sum_g w_{\rm exc}(g)&\le2n^2u,\label{eq:exc}\\
 \max_g w_{\rm reg}(g)&\le2u,\qquad
 \sum_g w_{\rm reg}(g)^2\le2n^2E_2(Q),\label{eq:reg}\\
 \sum_g w(g)&=n^2u^2.\label{eq:mass}
\end{align}
\end{lemma}
\begin{proof}
Consider first the additive setting. Write $f(x)=c+\lambda/(x+d)$, where $\lambda\ne0$. Translating the parameter sets absorbs $c$ and $d$ without changing $n$, $u$, or $E_2(Q)$; thus we may work with $h_{a,z}(x)=a+\lambda/(x+z)$. The pole and the value at infinity determine $z$ and $a$, respectively. For maps indexed by $(a,z)$ and $(b,r)$, put $t=z-r$. When $t=0$, their composition is the translation $x\mapsto x+a-b$, and such pairs have total weight $n^2u$. When $t\ne0$, the composition is
\begin{equation}\label{eq:addcomp}
 a+\lambda/t-\frac{\lambda^2/t^2}{x-b+\lambda/t}.
\end{equation}
The reciprocal coefficient, $-\lambda^2/t^2$, determines $t$ up to two possibilities in any characteristic (and uniquely when squaring is injective). Once $t$ is fixed, the pole $b-\lambda/t$ and the value at infinity $a+\lambda/t$ determine $b$ and $a$. For each triple $(a,b,t)$, exactly $r_Q(t)$ pairs $(z,r)\in Q^2$ satisfy $z-r=t$. Consequently, a regular transformation receives contributions from at most two triples and has weight at most $2u$. Summing the inequality $(v_1+v_2)^2\le2(v_1^2+v_2^2)$ over these triples yields
\[
 \sum_g w_{\rm reg}(g)^2
 \le 2\sum_{a,b}\sum_{t\ne0}r_Q(t)^2
 \le 2n^2E_2(Q).
\]
The argument remains valid in characteristic two.

In the multiplicative setting, write $f(x)=c f_0(dx)$ with $c,d\ne0$ and $f_0(x)=(x+1)/(x+\kappa)$, where $\kappa\ne0,1$. The scalings can be absorbed into $U$ and $Q$ without changing their multiplicative energies. The maps $a f_0(zx)$ are pairwise distinct: their value at infinity determines $a$, and their pole then determines $z$. With $t=z/r$, the composition $h_{a,z}h_{b,r}^{-1}$ has matrix
\begin{equation}\label{eq:matrix}
 G(a,b,t)=\begin{pmatrix}
 a(1-\kappa t)&ab(t-1)\\
 \kappa(1-t)&b(t-\kappa)
 \end{pmatrix}.
\end{equation}
For $t=1$, the composition is a dilation. If $t\ne1$, both off-diagonal entries are nonzero, and the projective invariant satisfies
\begin{equation}\label{eq:J}
 J(G)-1=\frac{G_{11}G_{22}}{G_{12}G_{21}}-1
 =-\frac{(\kappa-1)^2t}{\kappa(t-1)^2}.
\end{equation}
The right-hand side is nonzero and determines $t$ up to at most two possibilities, since the equation $t/(t-1)^2=v$ is quadratic and nontrivial in every characteristic, including characteristic two. Fix one such value of $t$. Since $\kappa(1-t)\ne0$, we normalize the matrix so that its lower-left entry equals $\kappa(1-t)$. The upper-right entry $ab(t-1)$ then determines $ab$. If $1-\kappa t\ne0$, we recover $a$ from $a=G_{11}/(1-\kappa t)$; if $t-\kappa\ne0$, we recover $b$ from $b=G_{22}/(t-\kappa)$. Together with $ab$, either value determines the other parameter. The only exceptional case is $1-\kappa t=t-\kappa=0$, equivalently $t=\kappa=\kappa^{-1}\ne1$. In characteristic different from two, this means $t=\kappa=-1$, for which the matrix is anti-diagonal and the map is of the form $x\mapsto c'/x$.

We therefore designate $t=1$ as exceptional, together with the additional fibre $t=\kappa=-1$ when it occurs. For fixed nonzero $t$, there are $r_Q(t)$ choices of $(z,r)$ and $n^2$ choices of $(a,b)$, so each exceptional fibre has total mass at most $n^2u$. On the regular fibres, each transformation corresponds to at most two triples $(a,b,t)$, each contributing multiplicity $r_Q(t)$. As in the additive case, this gives $w_{\rm reg}(g)\le2u$ and $\sum_g w_{\rm reg}(g)^2\le2n^2E_2(Q)$. Finally, \eqref{eq:mass} counts all ordered pairs of transformations.
\end{proof}

\begin{lemma}\label{lem:weighted}
With the notation of Lemma~\ref{lem:weights},
\begin{equation}\label{eq:weightedenergy}
 \sum_g w(g)i(g)\ll Mn^2u+M^{3/2}n^{3/2}uE_2(Q)^{1/4}
 +(1+M^2/p)n^2u^2.
\end{equation}
Over prime fields, the same bound holds with
$\sum_gw(g)(i(g)-\mu)$ on the left and $1+M^2/p$ replaced by $1$
on the right.
\end{lemma}
\begin{proof}
The exceptional contribution is at most $M\sum w_{\rm exc}$. For the regular contribution, decompose the weights dyadically by setting $\mathcal G_\tau=\{g:\tau\le w_{\rm reg}(g)<2\tau\}$, where $\tau$ ranges over powers of two. Lemma~\ref{lem:rich} controls the linear term by the total weight, leaving at most $M^{3/2}\sum_\tau \tau\card{\mathcal G_\tau}^{3/4}$. Set $W=n^2u^2$ and $V=2n^2E_2(Q)$. Then
\[
 \tau\card{\mathcal G_\tau}^{3/4}\le
 \min(W^{3/4}\tau^{1/4},V^{3/4}\tau^{-1/2}).
\]
Choose the splitting point $T=V/W=2E_2(Q)/u^2$. The elementary bounds $u^2\le E_2(Q)\le u^3$ give $2\le T\le2u$. Summing the first estimate over dyadic $\tau\le T$ and the second over $\tau\ge T$ yields, respectively, $O(W^{3/4}T^{1/4})$ and $O(V^{3/4}T^{-1/2})$. Both are $O(W^{1/2}V^{1/4})$, and $W^{1/2}V^{1/4}\ll n^{3/2}uE_2(Q)^{1/4}$, which proves \eqref{eq:weightedenergy}. The terms linear in $\card{\mathcal G_\tau}$ sum to at most $(1+M^2/p)\sum_gw(g)$. For the centered estimates, apply \eqref{eq:centersum}; negative summands may be discarded in an upper bound, while the residual linear error is $O(\sum_gw(g))$.
\end{proof}

\section{Pointwise energy transfer and proof of the expansion theorem}

\begin{proposition}[Asymmetric pointwise energy transfer]\label{prop:asymmetric}
Let $X,Y\subset\Gamma$ be nonempty finite sets of sizes $n,m$, respectively.
Assume either that $\Gamma=(F,+)$, $f$ is a nonaffine
M\"obius map and $Y$ avoids its pole, or that
$\Gamma=(F^\times,\cdot)$, $\alpha\beta\gamma\delta\ne0$,
and $f$ is defined and nonzero on $Y$. In either case
$\alpha\delta-\beta\gamma\ne0$.
Put $T=X\op f(Y)$ and $L=\card T$. Then
\begin{equation}\label{eq:asymgeneral}
 E_2(Y)\ll \frac{L^2m^{4/3}}{n^{2/3}}+\frac{L^2m^2}{p}.
\end{equation}
For $F=\F_p$ the stronger estimate is
\begin{equation}\label{eq:asymprime}
 E_2(Y)\ll \frac{L^2m^{4/3}}{n^{2/3}}
                 +\frac{m^3L}{p+1}.
\end{equation}
The characteristic-dependent term in \eqref{eq:asymgeneral} is omitted in characteristic zero. All implicit constants are absolute.
\end{proposition}
\begin{proof}
Apply Lemmas~\ref{lem:weights} and~\ref{lem:weighted} with $U=X$, $Q=Y^{-\op}$, and $S=T$.
Let
\[
 R_Y(x)=\card{\{(b,c)\in Y^2:b\op c=x\}},\qquad
 H(x)=\sum_{a\in X}\sum_{b\in Y}
       \one_T(h_{a,b^{-\op}}(x)).
\]
Extend $R_Y$ by zero to $\PP(F)$. For each representation $x=b\op c$ and each $a\in X$, the transformation $h_{a,b^{-\op}}$ maps $x$ to $a\op f(c)\in T$. Because these transformations are distinct, the resulting incidences are distinct for each fixed $x$. Therefore
\begin{equation}\label{eq:asympointwise}
 H(x)\ge nR_Y(x)\qquad(x\in\PP(F)).
\end{equation}
Put $E=E_2(Y)$. The abelian-group energy identity gives
$\sum_x R_Y(x)^2=E$, while projective bijectivity gives
\[
 \sum_x H(x)^2=\sum_gw(g)i_T(g),\qquad
 i_T(g)=\card{\{t\in T:g(t)\in T\}}.
\]
The identities remain valid over infinite fields, since the preimage of $T$ under each of the finitely many transformations is finite. Combining them with \eqref{eq:asympointwise} and \eqref{eq:weightedenergy} gives
\begin{equation}\label{eq:asymmaster}
 n^2 E\ll Ln^2m+L^{3/2}n^{3/2}mE^{1/4}
              +(1+L^2/p)n^2m^2.
\end{equation}
If the term involving $E^{1/4}$ is dominant, then
\[
 E^{3/4}\ll L^{3/2}n^{-1/2}m,
 \qquad E\ll L^2m^{4/3}n^{-2/3}.
\]
The remaining terms $Lm$ and $m^2$, arising after division of \eqref{eq:asymmaster} by $n^2$, are absorbed by the same bound: because $L\ge\max(n,m)$, one has $Lm\le L^2m^{4/3}n^{-2/3}$ and $m^2\le L^2m^{4/3}n^{-2/3}$. This proves \eqref{eq:asymgeneral}.

For the prime-field refinement, put $q=p+1$ and center the incidence function by defining $H_0(x)=H(x)-nmL/q$. Since $\sum_xR_Y(x)=m^2$, the pointwise estimate \eqref{eq:asympointwise} gives
\[
 \sum_xR_Y(x)H_0(x)\ge nE-\frac{nm^3L}{q}.
\]
If $E<2m^3L/q$, the required estimate \eqref{eq:asymprime} is immediate. Otherwise the preceding left-hand side is at least $nE/2$, and Cauchy--Schwarz yields
\[
 \frac{n^2E}{4}\le\sum_{x\in\PP(\F_p)}H_0(x)^2
 =\sum_gw(g)(i_T(g)-L^2/q).
\]
The centered estimate from Lemma~\ref{lem:weighted} consequently gives
\[
 n^2 E\ll Ln^2m+L^{3/2}n^{3/2}mE^{1/4}+n^2m^2.
\]
As in the uncentered argument, each potentially dominant term leads to $E\ll L^2m^{4/3}n^{-2/3}$. Together with the initial alternative, this proves \eqref{eq:asymprime}.
\end{proof}

\begin{proof}[Proof of Theorem~\ref{thm:main}]
In the reciprocal-sumset case take $\Gamma=(F,+)$ and
$f(x)=\lambda/x$; in the M\"obius-product case take
$\Gamma=(F^\times,\cdot)$ and $f$ as in the theorem.
In each case write $B=f(A)$ and $M=\card S$, so that
$\card B=N$. Apply Proposition~\ref{prop:asymmetric}
first with $X=Y=A$.
The diagonal case $n=m=N$, $L=M$ yields
\begin{equation}\label{eq:pointenergy}
 E_2(A)\ll M^2N^{2/3}+M^2N^2/p.
\end{equation}
Over prime fields it yields
\begin{equation}\label{eq:pointenergyprime}
 E_2(A)\ll M^2N^{2/3}+N^3M/(p+1).
\end{equation}
The same estimates hold for $B$: interchange $A$ and $B$ and
replace $f$ by its inverse M\"obius map. The admissibility
hypotheses are preserved under this operation.

Applying Cauchy--Schwarz to the representation function of
$S=A\op B$, followed by the energy identity, gives
\begin{equation}\label{eq:mixed}
 \frac{N^4}{M}\le \sum_t r_A(t)r_B(t)\le \sqrt{E_2(A)E_2(B)}.
\end{equation}
Here $r_B(t)=r_B(t^{-\op})$, and hence the same identity applies to both additive and multiplicative collisions. Substituting \eqref{eq:pointenergy} for each input set yields
\[
 N^4/M\ll M^2N^{2/3}+M^2N^2/p.
\]
According to which term dominates, we obtain $M^3\gg N^{10/3}$ or $M^3\gg pN^2$, proving \eqref{eq:main}. Over prime fields, \eqref{eq:pointenergyprime} instead yields
\[
 N^4/M\ll M^2N^{2/3}+N^3M/(p+1).
\]
The second term gives $M^2\gg (p+1)N$, proving \eqref{eq:prime}.
\end{proof}

\section*{AI Assistance Statement}
The author used ChatGPT, with GPT-6.0 available through a ChatGPT Plus
subscription, as an auxiliary research and writing tool during parts of
the exploratory development of this work. The tool was used for
brainstorming possible proof strategies, checking intermediate
calculations and logical dependencies, organizing arguments, and
assisting with the exposition and preparation of the manuscript.

All mathematical statements, proofs, references, and final formulations
appearing in this paper were independently reviewed and selected by the
author. The author assumes full responsibility for the correctness,
originality, and presentation of the results.

\end{document}